\documentclass[12pt,a4paper]{article}
\usepackage[latin2]{inputenc}
\usepackage{graphicx}
\usepackage{psfrag}
\usepackage[T1]{fontenc}
\usepackage{amsfonts}
\usepackage{amssymb}
\usepackage{amsmath,amsxtra,amsthm}

\newtheorem{lemma}{Lemma}
\newtheorem*{lemma*}{Lemma}
\newtheorem{theorem}{Theorem}
\newtheorem{observation}{Observation}

\newcommand{\newton}[2]{\left(\begin{array}{c}#1\\#2\end{array}\right)}
\newcommand{\classone}{{\mathcal C}^1}
\newcommand{\setN}{{\mathbb N}}
\newcommand{\setZ}{{\mathbb Z}}
\newcommand{\setR}{{\mathbb R}}
\newcommand{\setC}{{\mathbb C}}

\renewcommand{\author}{Miko\l{}aj Zalewski}
\renewcommand{\title}{Computer-assisteed proof of a periodic solution in a non-linear
feedback DDE}

\begin{document}

\noindent
\vskip1cm
\begin{center}
{\LARGE \bf \title}
\\[0.7cm]
{\large \author}\\
Institute of Mathematics, Jagiellonian University\\
Reymonta 5, 30-059 Krak\'o{}w, Poland\\[0.2cm]
E-mail: Mikolaj.Zalewski@im.uj.edu.pl
\\[0.7cm]
\end{center}

\vspace{0.3cm}

\begin{abstract}
In this paper we rigorously prove the existance of a non-trivial periodic orbit for the non-linear DDE: $x'(t) = K \sin(x(t-1))$ for $K=1.6$. We show that the equations on the Fourier equations have a solution by computing the local Brower degree. This degree can be
computed by using a homotopy which validity can be checked by checking a finite number of inequalities. Checking these inequalities is done by a computer program.
\end{abstract}

\section{Introduction}

We will show that the equation:
\begin{equation}
x'(t) = -K \sin\left(x(t-1)\right)
\end{equation}
has a periodic solution for $K = 1.6$. By a solution we mean a $\classone$ function that satisfies the equation.

Numerical simulations suggests that for $\frac{\pi}{2}<K<5.1$
there is a attracting periodic orbit oscillating around $0$. An
informal argument for the existence of such an orbit also comes
from the analysis of eigenvalues of the linearization
\cite{fizycy}. We will prove rigorously that for $K=1.6$ this
orbit exists.

In the remainder of the introduction we briefly outline our
method. We will use the Fourier coefficients of the periodic
orbit. First we rescale the time to obtain a $2\pi$ period. We
make a substitution: $\tilde{x}(t) = x\left(\frac{t}{\tau}\right)$
where $\tau$ is a parameter. If $\tau$ is equal to
$\frac{2\pi}{T}$ (where T is the period of the solution in the
original equation) we will have a $2\pi$-periodic solution of the
new equation:

\begin{align*}
\tilde{x}'(t) = \frac{1}{\tau}x'\left(\frac{t}{\tau}\right) =
-\frac{K}{\tau} \sin\left(x\left(\frac{t}{\tau}-1\right)\right) = \\
= -\frac{K}{\tau} \sin\left(x\left(\frac{t - \tau}{\tau}\right)\right)
= -\frac{K}{\tau} \sin(\tilde{x}(t - \tau))
\end{align*}

By renaming $\tilde{x}$ to $x$ we obtain:

\begin{equation}
\label{eq:rownanie}
x'(t) = -\frac{K}{\tau} \sin(x(t - \tau))
\end{equation}

We will prove for this equation that there exist a $\tau$ from a
small interval such there exist a $2\pi$-periodic orbit in a small
neighborhood of a specified function. Note that we don't obtain
the exact value of the period (even if numerical solution suggests
it is $4$, i.e. $\tau = \frac{\pi}{2}$) but treat $\tau$ as a
variable.

To prove the existence of the orbit we will use the method of
self-consistent bounds that was introduced in the context of
Kuramoto-Shivashinsky PDEs in \cite{zgliczyn2001},
\cite{zgliczyn}. When applied to the boundary value
problem for ODEs or DDEs this method is similar to the Cesari
method introduced in \cite{cesari} but doesn't require one of the
conditions -- see Section 2.4 in \cite{zgliczyn2001} for a
comparison. Below we briefly describe the method.

We will derive from (\ref{eq:rownanie}) the equation on Fourier coefficients.
Let's first write (\ref{eq:rownanie}) with a Taylor expansion instead of sinus:

\[ x'(t) = -\frac{K}{\tau} \sum_{n=0}^\infty \frac{(-1)^k}{(2k+1)!}[x(t-\tau)]^{2k+1}\]

The equations on the Fourier coefficients will look similar:

\[ \forall n \in \setZ: inc_n = -\frac{K}{\tau}\sum_{k=0}^{\infty}\frac{(-1)^k}{(2k+1)!}
\left(c^{*(2k+1)}\right)_n e^{-in\tau} \]
where $(c^{*(2k+1)})_n$ means convoluting the sequence $c$ with itself $2k+1$
times and then taking the $n$-th coefficient. The operation of convolution and
this notation is introduced in details in section \ref{sec:fourier}.

We will show that such a function has a non-trivial zero:

\[ F(\tau, c) = \left\{ in\tau e^{in\tau} c_n +
K\sum_{k=0}^{\infty}\frac{(-1)^k}{(2k+1)!}
\left(c^{*(2k+1)}\right)_n \right\}_{n=-\infty}^{\infty} \]

Where the domain is:

\begin{equation*}
X_\beta := \left\{x: \setZ \rightarrow \setC \middle| \forall n \in \setZ: |x_n| \leq \frac{\beta}{\left(|n|+1\right)^2}\mbox{ and }x_n = \overline{x_{-n}}\right\}
\end{equation*}
for some $\beta > 0$. On $X_\beta$ we consider the product topology (which is equivalent
to component-wise convergence). It turns out that on such a domain the sums in
the definition of $F$ are convergent, $F$ is continuous and each sequence $c \in X_\beta$
corresponds to a real-valued continuous function. Moreover if $F(c) = 0$ then the
function is $\classone$.

We will search for a zero in the neighborhood of an approximated
solution obtained from numerical simulations. We will denote by
$\hat{c}$ this approximation. For $n \leq n \leq 5$ it is equal
to:

\begin{equation*}
\begin{array}{|c|c|}
\hline
 n & \hat{c}_n \\
\hline
0 & 0 \\
1 & -0.1521000000-0.1163508047i \\
2 & 0 \\
3 & 0.0001123121-0.0002746107i \\
4 & 0 \\
5 & -0.0000008173-0.0000001014i \\
\hline
\end{array}
\end{equation*}
For $n>5$ the $\hat{c}_n$ is zero, for $n<0$ we have $\hat{c}_n = \overline{\hat{c}_{-n}}$.
The $\tau$ in the approximation is: $\hat{\tau} = 1.570796$.

For each $l \in \setN$ let us define the Galerkin projection $P_l$ and immersion $Q_l$:

\begin{equation*}
P_l: \setC^\setZ \ni c \rightarrow (c_{0}, \dots, c_l) \in \setR \times \setC^l
\end{equation*}
\begin{equation*}
Q_l: \setR \times \setC^{l} \ni (c_0, \dots, c_l)
\rightarrow (\dots, 0, 0, \overline{c_{-l}}, \dots, c_0, \dots, c_l, 0, 0, \dots)
\in \setC^\setZ
\end{equation*}

Let's note that we need only the non-negative terms in the finite space as the
negative terms can be obtained by conjugation. Also as $c_0 = \overline{c_{-0}}$
so we have $c_0 \in \setR$. From the compactness of $X_\beta$ it will be easy to show that:

\begin{lemma*} \textbf{\ref{lem:granica}}
Let $\beta>0$, $l_0 > 0$, $\underline{\tau}, \overline{\tau} \in \setR$ be fixed.
If for each $l > l_0$ there is a $c^l \in X_\beta$ and $\tau_l \in [\underline{\tau}, \overline{\tau}]$
such that $P_l F(\tau_l, c^l) = 0$ then there exist a $(c^0, \tau_0) \in X_\beta \times
[\underline{\tau}, \overline{\tau}]$ such that $F(\tau_0, c^0) = 0$.

Moreover if all the $c^l$ are in a closed set $D$ then $c^0 \in D$.
\end{lemma*}

Thus it's enough that there is a zero for every Galerkin projection of $F$:

\[ \tilde{F_l}: \setR \times \left(\setR \times \setC^l\right) \ni (\tau, c) \rightarrow
P_l F(\tau, Q_l(c)) \in \setR \times \setC^l \]

This will allow us to use a finite-dimensional topological method
to show a zero of $F$. Let's note that the real dimension of the
domain is $2l+2$ while of the image is $2l+1$ so if there is a
zero then one can expect a 1-dimensional manifold of zeros. This
manifold can be easily identified -- if $x(\cdot)$ is a solution
then $x(\cdot + \phi)$ is also. This means that if
$\tilde{F_l}(\tau, c_0, c_1, \dots, c_l) = 0$ then
$\tilde{F_l}(\tau, c_0, e^{i\phi}c_1, \dots, e^{il\phi}c_l) = 0$
what can be easy checked. To use the topological method we want to
have the zero isolated so we will limit the domain -- we assume
that $c_1 \in \hat{c}_1 + \setR$. Of course finding a zero in a
limited domain implies a zero of the full system. By $F_l$ we will
denote the $\tilde{F_l}$ with the smaller domain.

The before-mentioned method to prove the existence of a zero of
$F_l$ is the local Brouwer degree (introduced e.g. in \cite{S}). Let's denote
by $\deg(F_l, U, x)$ the degree of $x$ on $U$. It is known that if
the degree is non-zero then $\exists y \in U: f(y) = x$. We will
use a neighborhood of $\hat{c}$ as $U$ and $x=0$.

To compute the degree we will use the homotopy invariance of the local Brouwer degree. As our homotopy $H_l$ we will use a linear deformation of $F_l$ into a function $G_l$ that contains the most important terms of $F_l$ (that is not strictly a linearization of $F_l$ but it is close to):

\[ H_l(h, \tau, c) = h F_l(\tau, c) + (1-h) G_l(\tau, c)\]

The degree of $G_l$ will be easy to show that it is non-zero. What we will need to show is that $0 \notin H_l([0;1]; \partial U)$. To show that it's enough to show that the terms in $G_l$ dominates the (mainly non-linear) terms that are in $F_l$ but not $G_l$ -- that $\left| G_l \right| \geq \left| F_l - G_l\right|$. That part of the proof is computer-assisted -- the proofs of the estimates are in this paper but computing the exact values and checking that the inequalities holds is done by a computer program using the CAPD package for rigorous interval arithmetics.

\subsection{The program}
This program is written in C++ and can be downloaded from\\ http://www.im.uj.edu.pl/MikolajZalewski/dl/delay-sin.tgz .
The rounding-mode changing code required by the interval arithmetic is system-dependent and has been checked to work on PCs (both 32-bit and 64-bit) on both Windows (compiled with cygwin) and Linux (compiled with gcc). It should also work on SPARC and Mac OS X although that has not been tested. Using other CPUs or compilers might require modifications to the rounding code.

\section{Fourier coefficients}
\label{sec:fourier}

We will use the following notation: if $\{x_n\}_{n=-\infty}^\infty$ and
$\{y_n\}_{n=-\infty}^\infty$ are sequences with complex values then by
$x*y$ we will denote the convolution of the sequences:

\begin{equation*}
x*y := \left\{\sum_{k=-\infty}^{\infty} x_k y_{n-k}\right\}_{n = -\infty}^{\infty}
\end{equation*}

It is well known that if the sum converges the operation of convolution
is associative.
We will also use the notation $(...)_n$ for the $n$-th coefficient of the sequence in brackets -- e.g. $\left(x*y\right)_n = \sum_{k=-\infty}^{\infty} x_k y_{n-k}$. Also:

\begin{equation*}
x^{*k} := \underbrace{x * \cdots * x}_{k\mbox{ times}}\mbox{ (for $k \geq 1$)}
\end{equation*}

\begin{equation*}
(x^{*0})_n := \left\{ \begin{array}{l}
           1\mbox{ if }n=0 \\
           0\mbox{ if }n \neq 0
          \end{array} \right.
\end{equation*}

Of course the sum in the definition of the convolution may be not converging so we will limit our attention to a domain where the sum will be always convergent. The following lemma holds:

\begin{lemma}
\label{lem:szac}
If $x$ and $y$ is such that $\exists \alpha \geq 2, \beta_1, \beta_2>0: \forall n: |x_n| \leq
\frac{\beta_1}{\left(|n|+1\right)^\alpha}, |y_n| \leq
\frac{\beta_2}{\left(|n|+1\right)^\alpha}$, then $(x*y)_n$ is convergent for each $n$
and $\left|(x*y)_n\right| \leq C \frac{\beta_1\beta_2}{\left(|n|+1\right)^\alpha}$ where
 $C=\frac{2(2\alpha+1)}{\alpha-1}$
\end{lemma}

{\bf Proof:} Let's assume $n \geq 0$ and let's estimate $\left|(x*y)_n\right|$:

\[ \left| \sum_{k=-\infty}^{\infty} x_k y_{n-k} \right| \leq
\left| \sum_{k=0}^{\infty} x_{-k} y_{n+k}\right| +
\left| \sum_{k=1}^{n-1} x_{k} y_{n-k} \right| +
\left| \sum_{k=0}^{\infty} x_{n+k} y_{-k}\right|\]

\[
\begin{array}{rcl}
\left| \sum_{k=1}^{n-1} x_{k} y_{n-k} \right| & \leq &
\beta_1\beta_2 \sum_{k=1}^{n-1} \frac{1}{(k+1)^\alpha} \frac{1}{(n-k+1)^\alpha} \\
& \leq &
2\beta_1\beta_2 \sum_{1 \leq k \leq \frac{n}{2}} \frac{1}{(k+1)^\alpha} \frac{1}{(n-k+1)^\alpha} \\
& \leq &
2\beta_1\beta_2 \sum_{1 \leq k \leq \frac{n}{2}} \frac{1}{(k+1)^\alpha} \frac{1}{((n/2)+1)^\alpha} \\
& \leq &
2\frac{\beta_1\beta_2 2^\alpha}{(n+2)^\alpha} \sum_{1 \leq k \leq \frac{n}{2}} \frac{1}{(k+1)^\alpha} \\
& \leq &
2\frac{\beta_1\beta_2 2^\alpha}{(n+2)^\alpha} \sum_{k \geq 1} \frac{1}{(k+1)^\alpha} \\
\end{array}
\]

\[
\begin{array}{rcl}
\left| \sum_{k=0}^{\infty} x_{-k} y_{n+k} \right| & \leq &
\beta_1\beta_2 \sum_{k=0}^{\infty} \frac{1}{(k+1)^\alpha} \frac{1}{(n+k)^\alpha+1} \\
& \leq &
\frac{\beta_1\beta_2}{(n+1)^\alpha} \sum_{k \geq 0} \frac{1}{(k+1)^\alpha} \\
\end{array}
\]

Analogically:
\[
\left| \sum_{k=0}^{\infty} x_{n+k} y_{-k} \right| \leq
\frac{\beta_1\beta_2}{(n+1)^\alpha} \sum_{k \geq 0} \frac{1}{(k+1)^\alpha}
\]

The sums $\sum_{k=k_0}^{\infty} \frac{1}{k^\alpha+1}$ can be estimated by integrals:

\[ \sum_{k=0}^{\infty} \frac{1}{(k+1)^\alpha}\leq 1 + \int_{0}^\infty \frac{1}{(x+1)^\alpha} dx
= \frac{\alpha}{\alpha-1}\]
\[ \sum_{k=1}^{\infty} \frac{1}{(k+1)^\alpha}\leq \frac{1}{2^\alpha} + \int_{1}^\infty \frac{1}{(x+1)^\alpha} dx
= \frac{\alpha+1}{2^\alpha(\alpha-1)}\]

From that we obtain:

\[ \left| \sum_{k=-\infty}^{\infty} x_k y_{n-k} \right| \leq
\frac{\beta_1\beta_2}{(n+1)^\alpha}\left[ 2^{\alpha+1} \frac{\alpha+1}{2^\alpha(\alpha-1)}
+ 2 \frac{\alpha}{\alpha-1} \right] \leq
\frac{\beta_1\beta_2}{(n+1)^\alpha} \frac{2(2\alpha+1)}{\alpha-1} \]

The result for $n<0$ can be obtained by analogous estimations or by taking a sequences
$\tilde{x}$, $\tilde{y}$: $\tilde{x}_n := x_{-n}$, $\tilde{y}_n := y_{-n}$ and applying for
them the result for $n>0$.

\qed

With one exception we will use this lemma for $\alpha=2$. For $\alpha=2$
we have $C=10$.

\begin{observation}
If $x$ is such that $\exists \alpha \geq 2, \beta>0: \forall n: |x_n| \leq
\frac{\beta}{\left(|n|+1\right)^\alpha}$ then $\left|(x^{*k})_n\right| \leq C^{k-1} \frac{\beta^k}{\left(|n|+1\right)^\alpha}$ where
 $C=\frac{2(2\alpha+1)}{\alpha-1}$
\end{observation}

\begin{lemma}
\label{lem:splotrzecz}
If $x$ and $y$ are such that $x_{-n} = \overline{x_n}$, $y_{-n} = \overline{y_n}$ then
$(x*y)_{-n} = \overline{(x*y)_n}$
\end{lemma}

Thus we have that if $x \in X_{\beta_1}$ and $y \in X_{\beta_2}$ then $x * y \in X_{C\beta_1\beta_2}$
(where $X_\beta$ was defined in the introduction and $C$ is from Lemma \ref{lem:szac}).
It will be also useful to define a set of sequences from any $X_\beta$:

\begin{equation*}
X:= \bigcup_{\beta > 0} X_\beta
\end{equation*}

We will limit ourselves to the sequences in $X$. For them we have:
\begin{lemma}
Let $c \in X$. Then $\sum_{n=-\infty}^\infty c_n e^{int}$ converges to a real-valued continuous function
\end{lemma}

{\bf Proof: } The functions $\sum_{n=-N}^N c_n e^{int}$ are continuous and real-valued as $c_n = \overline{c_{-n}} \Rightarrow c_ne^{int} + c_{-n}e^{-int} \in \setR$. They uniformly convergent as $\sum_{n=-\infty}^\infty \left|c_n e^{int}\right| \leq \sum_{n=-\infty}^\infty \frac{\beta}{(|n|+1)^2} < \infty$. Thus the limit is also real-valued and continuous. \qed

Let's note that $c \in X$ doesn't guarantee that the function is $\classone$:

We use the convolution because of this well known fact:
\begin{lemma}
If $c \in X$ are the Fourier coefficients of $x(t)$, $d \in X$ are
the coefficients of $y(t)$ then the Fourier coefficients of
$x(t)\cdot y(t)$ are $c*d$. As a consequence the coefficients of
$x^n(t)$ are $c^{*n}$.
\end{lemma}

Now we can write the equation (\ref{eq:rownanie}) on the Fourier
coefficients.

\begin{theorem}
\label{tw:fourier}
Let $\tau$ be fixed and $x(t): \setR \rightarrow \setR$ be a $2\pi$-periodic function
with the Fourier coefficients $c \in X$. Then:
\begin{itemize}
 \item[(i)] for each $n \in \setZ$ the sum on the right-hand side of the following equation converges:

\begin{equation}
\label{eq:nawsp}
inc_n = -\frac{K}{\tau}\sum_{k=0}^{\infty}\frac{(-1)^k}{(2k+1)!} \left(c^{*(2k+1)}\right)_n e^{-in\tau}
\end{equation}

\item[(ii)] $x(t)$ is a $2\pi$-periodic solution of (\ref{eq:rownanie}) $\Leftrightarrow$ $c$
satisfies the equations ($\ref{eq:nawsp}$).

\end{itemize}

\end{theorem}

To prove the theorem we will need some lemmas. As it has been already mentioned
that $c \in X$ doesn't imply that $x(t)$ is $\classone$. However if $c$ is the solution
of equation (\ref{eq:nawsp}) then we have the following lemma that will allow us to
show that $x(t)$ is $\classone$ ($x(t)$ is even ${\cal C}^\infty$):

\begin{lemma}
\label{lem:zwiwykl}
If $c$ satisfies equation (\ref{eq:nawsp}) and $\beta > 0, \alpha \geq 2$ are such that
$\forall n: \left| c_n \right| \leq \frac{\beta}{(|n|+1)^\alpha}$ then
$\exists \beta': \forall n: \left| c_n \right| \leq \frac{\beta'}{(|n|+1)^{\alpha+1}}$
\end{lemma}

{\bf Proof: }
Let $C := \frac{2(2\alpha+1)}{\alpha-1}$. We have that:

\begin{equation}
\label{eq:szaczb}
 \begin{array}{rcl}
\left| -\frac{K}{\tau}\sum_{k=0}^{\infty}\frac{(-1)^k}{(2k+1)!} \left(c^{*(2k+1)}\right)_n e^{-in\tau} \right| & \leq & \frac{K}{\tau}
\sum_{k=0}^{\infty} \frac{1}{(2k+1)!} | \left(c^{*(2k+1)}\right)_n | \\
& \leq & \frac{K}{\tau} \sum_{k=0}^{\infty} \frac{1}{(2k+1)!} C^{2k}\frac{\beta^{2k+1}}{(|n|+1)^\alpha} \\
& = & \frac{K}{\tau} \frac{1}{C(|n|+1)^\alpha} \sinh (C\beta)
   \end{array}
\end{equation}

We also have $\left|inc_n\right|=|n| \left|c_n\right|$. The two sides of the equation must be equal hence we obtain:
\[ \frac{K}{\tau} \frac{1}{C(|n|+1)^\alpha} \sinh (C\beta) \geq |n| \left|c_n\right| \]
Thus for $n \neq 0$ the $\beta' = \frac{2 K}{C\tau} \sinh (C\beta)$ satisfies the assertion. If
it is not satisfied for $n=0$ we can increase the $\beta'$.
\qed

Increasing the $\alpha$ is important as we have:
\begin{lemma}
\label{lem:fourierc1}
Let $c$ be a sequence of complex values satisfying $\left| c_n \right| \leq \frac{\beta}{(|n|+1)^3}$ for some $\beta$ and $c_n = \overline{c_{-n}}$. Then the sequence $c$ is a sequence of Fourier coefficients of a real-valued $\classone$ function $x(t)$.
\end{lemma}

{\bf Proof:}
The sequences $\sum_{k=-n}^n c_k e^{ikt}$ and $\sum_{k=-n}^n ikc_k e^{ikt}$ are real-valued, the second is the derivative of the first one and are uniformly convergent as $n \rightarrow \infty$. Hence both converge to continuous function and the second is the derivative of the first one. Thus $x(t) = \sum_{k=-\infty}^\infty c_k e^{ikt}$ is $\classone$.
\qed

{\bf Proof of the Theorem \ref{tw:fourier}: }

{\it Ad (i): } The convolutions are convergent because $x \in X$.
From the equation (\ref{eq:szaczb}) from the proof of Lemma
\ref{lem:zwiwykl} we have that if $x \in X_\beta$ then
$\sum_{k=0}^{\infty}\left|\frac{(-1)^k}{(2k+1)!}
\left(c^{*(2k+1)}\right)_n e^{-in\tau} \right| \leq
\frac{1}{C(|n|+1)^2} \sinh (C\beta) < \infty$

{\it Ad (ii): } Implication $\Rightarrow$: It's enough to show that $\left\{inc_n\right\}_{n=-\infty}^\infty$ are the Fourier coefficients of $x'(t)$ while
$-\frac{K}{\tau}\sum_{k=0}^{\infty}\frac{(-1)^k}{(2k+1)!} \left(c^{*(2k+1)}\right)_n e^{-in\tau}$ are the Fourier coefficients of $-\frac{K}{\tau} \sin x(t-\tau)$.

The first part can be obtained by integrating $\int_0^{2\pi} x'(t) e^{ikt} dt$ by parts.

As for the second we have that $-\frac{K}{\tau} \sum_{k=0}^N \frac{(-1)^k}{(2k+1)!} \left(x(t-\tau)\right)^{2k+1} \rightarrow -\frac{K}{\tau} \sin x(t-\tau)$ as $N \rightarrow \infty$. The Fourier coefficients of $x(t-\tau)$ are equal to $d:=\left\{c_n e^{-in\tau}\right\}_{n=-\infty}^\infty$. Hence the coefficients of $x(t-\tau)^{2k+1}$ are $d^{*(2k+1)}$ which is equal to $c^{*(2k+1)}e^{-in\tau}$.

Hence we have that the Fourier coefficients of $-\frac{K}{\tau} \sum_{k=0}^N \frac{(-1)^k}{(2k+1)!} \left(x(t-\tau)\right)^{2k+1}$ are equal to
$\left\{-\frac{K}{\tau}\sum_{k=0}^{N}\frac{(-1)^k}{(2k+1)!} \left(c^{*(2k+1)}\right)_n e^{-in\tau}\right\}_{n=-\infty}^\infty$.
It has been shown that this sequence in convergent as $N \rightarrow \infty$ thus the $n$-th coefficient of $-\frac{K}{\tau} \sin x(t-\tau)$ is $-\frac{K}{\tau}\sum_{k=0}^{\infty}\frac{(-1)^k}{(2k+1)!} \left(c^{*(2k+1)}\right)_n e^{-in\tau}$. This ends the proof of this case.

Implication $\Leftarrow$: from Lemmas \ref{lem:zwiwykl} and
\ref{lem:fourierc1} we obtain that $x(t)$ is a $\classone$
function. We know that
$-\frac{K}{\tau}\sum_{k=0}^{\infty}\frac{(-1)^k}{(2k+1)!}
\left(c^{*(2k+1)}\right)_n e^{-in\tau}$ are the Fourier
coefficients of $-\frac{K}{\tau} \sin x(t-\tau)$. They are equal
to $inc_n$ -- the Fourier coefficients of $x'(t)$. So both
$-\frac{K}{\tau} \sin x(t-\tau)$ and $x'(t)$ are continuous and
$2\pi$-periodic functions with equal Fourier coefficients. Hence
the functions themselves are equal and equation
(\ref{eq:rownanie}) is satisfied. \qed

Let's define a function $F: \setR \times X \rightarrow \setC^\setZ$ (with the product topology on $\setC^\setZ$) corresponding to the equation ($\ref{eq:nawsp}$):

\begin{equation*}
F(\tau, c) := \left\{in\tau e^{in\tau}c_n + K\sum_{k=0}^{\infty} \frac{(-1)^k}{(2k+1)!} \left(c^{*(2k+1)}\right)\right\}_{n=-\infty}^\infty
\end{equation*}

Of course having $F(\tau, c) = 0$ is equivalent to the fact that $\tau, c$
satisfies equation (\ref{eq:nawsp}). In the rest of the paper we will show that $F$ has
a nontrivial zero. Let's note that $F(\tau, c)_n = \overline{F(\tau, c)_{-n}}$.

As we will use topological tools we will want $F$ to be
continuous. First let's note that the convolution is not
continuous on the whole $X$ (as written above we use the product
topology on $X$) but we have:

\begin{lemma}
The operation of convolution is continuous on each $X_\beta$
\end{lemma}

{\bf Proof: } Let's fix some $x^0, y^0 \in X_\beta$. Let's fix
some $\delta$ and take some $x, y \in X_\beta$ from a neighborhood
of $\left(x^0, y^0\right)$: $x, y$ are such that $\forall |n| \leq
N_\delta: |x_n- (x^0)_n|, |y_n - (y^0)_n| < \delta$ where
$N_\delta$ is large enough that for $|n| > N_\delta: x, y, x^0,
y^0 \in X_\beta \Rightarrow |x_n - (x^0)_n|, |y_n - (y^0)_n| <
\delta$. We have that:

\[ \left|\left(x^0*y^0 - x * y\right)_n\right| \leq
\left| \left(x^0*(y^0-y)\right)_n\right| + |\left((x^0-x)*y\right)_n|\]

\[ \left|\left(x^0*(y^0-y)\right)_n\right| \leq \sum_{j=-\infty}^\infty \left| \left(x^0\right)_j\right|\cdot \left|\left(y^0\right)_{n-j} - y_{n-j}\right| \leq \beta\delta \sum_{j=-\infty}^\infty\frac{1}{(|j|+1)^2}\]
This tends to zero as $\delta \rightarrow 0$.
After an analogous estimation for $|\left((x^0-x)*y\right)_n|$ we have that
the convolution is continuous.
\qed

\begin{lemma}
The function $F$ is continuous on each $X_\beta$.
\end{lemma}

{\bf Proof:} The convolutions are continuous on each $X_\beta$ and
the result of a convolution lays in a $X_{\beta'}$ for some
$\beta'$ so $K\sum_{k=0}^{N}\frac{(-1)^k}{(2k+1)!}
\left(c^{*(2k+1)}\right)_n$ is continuous for each $N$. From
estimations as in equation (\ref{eq:szaczb}) we have that for each
coefficient this series converges uniformly as $N \rightarrow
\infty$ hence we obtain that the limit
$K\sum_{k=0}^{\infty}\frac{(-1)^k}{(2k+1)!}
\left(c^{*(2k+1)}\right)_n$ is continuous.

The term $in\tau e^{in\tau}c_n$ is also continuous so $F$ is continuous.

\qed

\section{Some estimates}

As mentioned in the introduction we will need some estimates to
show that the inequality holds in the neighborhood of $\hat{c}$.
We will use two kinds of sets: $\hat{c}+X_{\beta_2}$ that will be
used to obtain finer estimates and $X_{\beta_1}$ (where $\beta_1$
will be large enough to contain the whole set
$\hat{c}+X_{\beta_2}$) for some more rough but simpler ones.

In this section we will only assume of $\hat{c}$ that almost all coefficients are equal to zero. By $Y_l$ we will denote the space of the possible values of $\hat{c}$ -- the set of sequences such that at most the elements $-l, \dots, l$ are non-zero:

\[ Y_l := \{c \in X: \forall n: |n| > l \Rightarrow c_n = 0\} \]

First let's note two simple properties:

\begin{lemma}
If $c \in Y_l$ then $c^{*k} \in Y_{kl}$
\end{lemma}

\begin{lemma}
\label{lem:splsumy}
If $x, y \in X$ then
\begin{equation*}
(x+y)^{*n} = \sum_{k=0}^{k=n} \newton{n}{k} x^{*k} * y^{*(n-k)}
\end{equation*}
\end{lemma}

To estimate $\frac{K}{p!} \left(x^{*p}\right)_n$ -- an element of the sum in the equation ($\ref{eq:nawsp}$) -- for $x \in X_\beta$ it's enough to use Lemma \ref{lem:szac}. However for the sets of the form $c+X_\beta$ we will use a more sophisticated estimation:
\begin{lemma}
\label{lem:srodek}
Let $c \in Y_l$, $p>1$, $\beta>0$. Then for any $x \in X_\beta$:

\begin{equation*}
\left|\frac{K}{p!} (c+x)^{*p}_n\right| \leq \frac{K}{p!} \left( \left|c^{*p}_n\right| +
\sum_{k=0}^{p-1} \newton{p}{k} \sum_{j=-lk}^{lk} \left|c^{*k}_j\right| C^{p-k-1} \frac{\beta^{p-k}}{(|n-j|+1)^2}\right)
\end{equation*}
Where $C=10$

\end{lemma}

{\bf Proof:} From Lemmas \ref{lem:szac} and \ref{lem:splsumy} we have:

\begin{equation*}
\begin{array}{rcl}
\left|\frac{K}{p!} (x+c)^{*p}_n\right| & \leq & \frac{K}{p!} \sum_{k=0}^{p}
\newton{p}{k} \left|\left(c^{*k} * x^{*(p-k)}\right)_n\right| \\
& \leq & \frac{K}{p!} \sum_{k=0}^{p} \newton{p}{k}
\sum_{j=-\infty}^{\infty} \left|c^{*k}_j\right| \left|x^{*(p-k)}_{n-j}\right| \\
& \leq & \frac{K}{p!} \left( \left|c^{*p}_n\right| +
\sum_{k=0}^{p-1} \newton{p}{k} \sum_{j=-\infty}^{\infty} \left|c^{*k}_j\right| C^{p-k-1} \frac{\beta^{p-k}}{(|n-j|+1)^2}\right)
\end{array}
\end{equation*}

We have $c \in Y_l$ so $c^{*k}_j = 0$ for $j>kl$ what ends
the proof.
\qed

Of course we can use the previous lemma only for a finite number
of terms. To estimate the tail we will use a weaker estimation by
using a neighborhood of the second type, applying Lemma
\ref{lem:szac} and summing the geometric sequence:
\begin{lemma}
\label{lem:ogon}
Let $N$ be odd, $\beta \in [0; \frac{N}{10}]$. If $x \in X_\beta$  then
\begin{equation*}
\left|\left(\sum_{k=\frac{N-1}{2}}^{\infty}\frac{K}{(2k+1)!}x^{*(2k+1)}\right)_n\right|
\leq
\frac{K \cdot (C\beta)^{N-1}}{N!(1-\left(\frac{C\beta}{N}\right)^2)} \cdot \frac{\beta}{\left(|n|+1\right)^2}
\end{equation*}
Where $C = 10$.
\end{lemma}

{\bf Proof:} Using Lemma {\ref{lem:szac}}:

\begin{equation*}
\begin{array}{rcl}
\left|\left(\sum_{k=\frac{N-1}{2}}^{\infty}\frac{K}{(2k+1)!}x^{*(2k+1)}\right)_n\right|
 & \leq & \sum_{k=(N-1)/2}^{\infty}
     \frac{K}{N!N^{2k-(N-1)}}\left|\left(x^{*(2k+1)}\right)_n\right| \\
 & \leq & \frac{K}{N!} \sum_{k=(N-1)/2}^{\infty}
            \frac{(C\beta)^{2k}}{N^{2k-(N-1)}} \cdot \frac{\beta}{\left(|n|+1\right)^2} \\
 & = & \frac{K(C\beta)^{N-1}}{N!} \left(\sum_{k=0}^{\infty} \left(\frac{C^2\beta^2}{N^2}\right)^k\right)
                  \cdot \frac{\beta}{\left(|n|+1\right)^2} \\
 & = & \frac{K(C\beta)^{N-1}}{N!\left(1-\left(\frac{C\beta}{N}\right)^2\right)}
                  \cdot \frac{\beta}{\left(|n|+1\right)^2}
\end{array}
\end{equation*}

The geometric sequence is convergent because $\beta < \frac{N}{10} \Rightarrow
\left(\frac{C\beta}{N}\right)^2 < 1$.
\qed

From Lemma \ref{lem:srodek} we can obtain an estimation that after
multiplication by $(|n|+1)^2$ is independent of $n$. That will be
used to estimate the term in all inequalities for high $n$-s with
one formula:
\begin{lemma}
\label{lem:wysoko}
Let $\beta > 0$, $c \in Y_l$, $p>1$,
$N > pl$. Then for each $x \in X_\beta$, $n > N$ we have:
\begin{equation*}
\left|\frac{K}{p!} (x+c)^{(*p)}_n\right| \leq \frac{K}{p!} \left(
\sum_{k=0}^{p-1} \newton{p}{k} \sum_{j=-lk}^{lk} \left|c^{*k}_j\right| (C\beta)^{p-k-1} \frac{(N+1)^2}{(N-j+1)^2}\right)\frac{\beta}{(|n|+1)^2}
\end{equation*}

\end{lemma}

{\bf Proof:} Note that $n > pl \Rightarrow c^{*p}_n = 0$ and
$\frac{\beta}{(|n-j|+1)^2} = \frac{(n+1)^2}{(|n-j|+1)^2} \cdot \frac{\beta}{(n+1)^2} \leq \frac{(N+1)^2}{(N-j+1)^2} \cdot \frac{\beta}{(n+1)^2}$
as $N>j$ and apply Lemma \ref{lem:srodek}
\qed

For the first terms we will need to have a better estimation than in Lemma \ref{lem:srodek}
so we will regroup the terms (to understand why this regrouping helps let's
compare two estimates: $|1.1\cdot x -1 \cdot x| \leq 1.1|x|+1|x| = 2.1|x|$ and
$|1.1\cdot x - 1\cdot x| = |(1.1-1) x| = 0.1|x|$).

\begin{lemma}
\label{lem:grupowanie}
Let $x \in X_\beta$, $c \in Y_l$.

\begin{equation*}
\begin{array}{ll}
in\tau e^{in\tau}(c_n+x_n) + K\sum_{k=0}^3 &\frac{(-1)^k}{(2k+1)!}(c+x)^{*(2k+1)}_n = \\
& \left(in\tau e^{in\tau}+K\right)c_n
- \frac{K}{3!}\left(c^{*3}\right)_n
+ \frac{K}{5!}\left(c^{*5}\right)_n
- \frac{K}{7!}\left(c^{*7}\right)_n +\\
&+ \left(in\tau e^{in\tau}+K\right)x_n +
\sum_{p=1}^{7}\sum_{j=-6l}^{6l} \gamma_{p, j} x^{*p}_{n-j}
\end{array}
\end{equation*}
Where $\gamma_{p, j} = -\frac{K}{3!}\newton{3}{p}\left(c^{*(3-p)}\right)_j
+\frac{K}{5!}\newton{5}{p}\left(c^{*(5-p)}\right)_j
-\frac{K}{7!}\newton{7}{p}\left(c^{*(7-p)}\right)_j$ and we assume $\newton{n}{k}=0$
for $k>n$.
\end{lemma}

\section{Proving the existence of the periodic solution}
\label{sec:prove}

We will search for the orbit in the neighborhood of $\hat{c}$ and
$\hat{\tau}$ defined in the introduction. Let's define the sets
and boundaries on which we will operate:
\[ X_1 := X_{\beta_1} \]
\[ X_2 := \hat{c}+ X_{\beta_2} \]
\[ X_3 := \{y \in X_2: y_1-\hat{c}_1 \in \setR\} \]
\[ \underline{\tau}:= \hat{\tau} - \Delta\tau \]
\[ \overline{\tau} := \hat{\tau} + \Delta\tau \]
where $\Delta\tau = 0.000001$, $\beta_2 = 0.0000002438$, $\beta_1 = 0.766763$.

We will prove that a solution exists in the set $[\underline{\tau}; \overline{\tau}] \times X_3$. The $\beta_1$, $\beta_2$ are such that $X_3 \subset X_2 \subset X_1$.

As mentioned in the introduction we will use the Galerkin projections of $F$ with the condition $y_1 - \hat{c}_1 \in \setR$ added to make the zero isolated. The $P_l$ and $Q_l$ were defined in the introduction as the projection and immersion of the finite-dimensional space. Let's first prove the lemma stated in the introduction.

\begin{lemma}
\label{lem:granica}
Let $\beta>0$, $l_0>0$, $\underline{\tau}, \overline{\tau} \in \setR$ be fixed.
If for each $l>l_0$ there is a $c^l \in X_\beta$ and $\tau_l \in [\underline{\tau}, \overline{\tau}]$
such that $P_l F(\tau_l, c^l) = 0$ then there exist $(c^0, \tau_0) \in X_\beta \times
[\underline{\tau}, \overline{\tau}]$ such that $F(\tau_0, c^0) = 0$.

Moreover, if all the $c^l$ are in a closed set $D$ then $c^0 \in D$.
\end{lemma}

{\bf Proof:} From the compactness of $[\underline{\tau},
\overline{\tau}] \times X_\beta$ there exists a subsequence $l_k$
such that $c^{l_k}$  converges to a limit $c^0$ and $\tau_{l_k}$
converges to $\tau_0$. Let's fix $n \geq 0$ and note that for $l >
n$ if $P_l(F(\tau, c)) = 0$ then $F(\tau, c)_n = 0$, so from the
continuity of $F$ we have $F(\tau, c^0)_n = \lim_{k \rightarrow
\infty} F(\tau_{l_k}, c^{l_k})_n = \lim_{k \rightarrow \infty}
P_{l_k} F(\tau_{l_k}, c^{l_k})_n = 0$. For $n<0$ we have
$F(\tau_0, c^0)_n = \lim_{k \rightarrow \infty} F(\tau_{l_k},
c^{l_k})_n = \overline{\lim_{k \rightarrow \infty} F(\tau_{l_k},
c^{l_k})_{-n}} = 0$ Thus $F(\tau_0, c^0) = 0$.

The last assertion is obvious.
\qed

Before defining the homotopy let's introduce some auxiliary
notations. Let's denote: $f_n(\tau):=in\tau e^{in\tau}$. By $L_n$
we will denote the linear part of $f_n$: $L_n(\tau) :=
(ine^{in\tau_0}-n^2e^{in\tau_0})(\tau - \tau_0)$. By $r_n$ we will
denote the non-linear part: $r_n(\tau) := f_n(\tau) - f_n(\tau_0)
- L_n(\tau)$.

We will define the homotopy on the whole infinite-dimensional space: $H : [0;1] \times [\underline{\tau}; \overline{\tau}] \times X_3 \rightarrow \setC^\setZ$. We will use only the Galerkin projections of this  homotopy to deform $F_l$ (that is the Galerkin projection of $F$ -- see introduction) but to prove that there are no zeros on the boundaries it's a bit more convenient to define the homotopy on the whole space. As written in the introduction the homotopy is a linear deformation of $F$ to a nearly linear function $G$:

\[ H(h, \tau, x) := hF(\tau, \hat{c}+x) + (1-h)G(\tau, x) \]
where $G$ on the $n$-th coefficient is equal to:
\[
 G(\tau, x)_n := \left\{
    \begin{array}{l}
    \left(f_n(\tau)+K\right)x_n\mbox{ for }n \neq \pm 1 \\
    \left(f_n(\tau_0)+K+\gamma_{1,0}+\gamma_{1,2n}\right)x_n + L(\tau)\hat{c}_n\mbox{ for } n = \pm 1
    \end{array}\right.
\]
Where $\gamma_{p, j}$ is from Lemma \ref{lem:grupowanie}.

The $G$ is not strictly a linearization of $F$ as it doesn't contain all the linear
terms and is not linear with respect to $\tau$. However it contains the most important
terms -- the rest will be shown to be small compared to them. Let's note that for $n = \pm1$
we have $x_{n-2n}=\overline{x_n}=x_n$ (because $x_{\pm 1} \in \setR$) so the term
with $\gamma_{1,2n}$ in $F$ will be a linear term with respect to $x_n$.

One can write explicit formulas for $H$. For $n \neq \pm 1$ we have:
\begin{equation*}
\begin{array}{l}
 H(h, \tau, x)_n := \left(f_n(\tau)+K\right)x_n +\\
~+ h\left((f_n(\tau)+K)\hat{c}_n +
\sum_{k=1}^{\infty} \frac{K}{(2k+1)!} \left((x+\hat{c})^{*(2k+1)}\right)_n\right)
\end{array}
\end{equation*}

And for $n = \pm 1$:
\begin{equation*}
\begin{array}{rcl}
 H(t, \tau, x)_n & := & \left(f_n(\tau_0)+K+\gamma_{1, 0} +\gamma_{1,2n}\right)x_n
+ L_n(\tau)\hat{c}_n +  \\
 & & + h\left(\left(f_n(\tau_0)+K\right)\hat{c}_n
- \frac{K}{3!}\left(\hat{c}^{*3}\right)_n
+ \frac{K}{5!}\left(\hat{c}^{*5}\right)_n
- \frac{K}{7!}\left(\hat{c}^{*7}\right)_n\right) +\\
 & & + h\sum_{-6\cdot 5 \leq j \leq 6\cdot 5, j \neq 0,2n} \gamma_{1, j} x_{n-j} +
h\sum_{p=2}^{7}\sum_{j=-6\cdot 5}^{6\cdot 5} \gamma_{p, j} x^{*p}_{n-j} + \\
 & & + hr_n(\tau)\hat{c}_n + h(f_n(\tau) - f_n(\tau_0))x_n + h R(\hat{c}+x)
\end{array}
\end{equation*}
Where $R$ is a short notation for:
\[ R(y) := \sum_{k=4}^\infty \frac{K}{(2k+1)!} y^{*(2k+1)}\]

To show that any Galerkin projection of $H$ doesn't have a zero on the boundary
of $[\underline{\tau}; \overline{\tau}] \times P_l(X_3)$ it's enough to show that $H$ doesn't have a zero on the boundary. More exactly:

\begin{theorem}
Let $x \in X_3$, $\tau \in [\underline{\tau}, \overline{\tau}]$.
Let $x$ be such that $\exists n: |x_n| = \frac{\beta_2}{(|n|+1)^2}$ or
let $\tau \in \{\underline{\tau}, \overline{\tau}\}$.
Then $\forall h \in [0;1]: H(h, \tau, x) \neq 0$.
\end{theorem}

{\bf Proof:}
The proof is computer assisted -- the calculations are done by the program.
The computation is as follow.

First let's note that if for $n < 0$ we have $|x_n| =
\frac{\beta_2}{(|n|+1)^2}$ then this is also true for $-n$. Thus
it's enough to consider $n \geq 0$ (and $\tau$). We do it in two
steps:

{\bf 1.} Let's assume that $|x_n| = \frac{\beta_2}{(|n|+1)^2}$ for some $n \geq 0$, $n \neq 1$.
We will prove that:

\[ \begin{array}{rcl} \left|\left(f_n(\tau)+K\right)x_n\right| & > &
  \left|(f_n(\tau)+K)\hat{c}_n +
  \sum_{k=1}^{\infty} \frac{K}{(2k+1)!} \left((x+\hat{c})^{*(2k+1)}\right)_n\right| \end{array} \]

As the value of $H(h, \tau, x)_n$ is the difference of these two terms (with the second one multiplied by a $h<1$), it will follow that $H(h, \tau, x)_n \neq 0$ thus $H(h, \tau, x) \neq 0$.
The LHS will be estimated by:

\begin{equation}
\label{eq:homotoplhs}
 \left|\left(f_n(\tau)+K\right)x_n\right| \geq \left| \left|f_n(\tau)\right| - K\right|
 \frac{\beta_2}{(|n|+1)^2}
\end{equation}

and $\left|f_n(\tau)\right| = n|\tau|$ is an interval in the interval arithmetic.

The terms of RHS are estimated depending on $n$. For $n < 225$ we can estimate each $n$
separately. For $n \geq 225$ we will want to check all the inequalities in some finite computations.

Let's note that for $n \geq 255$ we have $(f_n(\tau)+K)\hat{c}_n = 0$ as
$\hat{c}_n = 0$. In the sum $\sum_{k=1}^{\infty} \frac{K}{(2k+1)!}
\left((x+\hat{c})^{*(2k+1)}\right)_n$ we will estimate the terms $k=1,
\dots, 22$ from Lemma \ref{lem:wysoko} while the term $k>22$ from
Lemma \ref{lem:ogon}. Now if we multiply the estimates for both
LHS and RHS by $(|n|+1)^2$ they will be independent of $n$. Hence
we will check all the $n \geq 225$ by checking only one
inequality.

For $n < 225$ we have to use some more sophisticated estimations.
Let's use Lemma \ref{lem:grupowanie} to group the terms:

\begin{equation*}
\begin{array}{rcl}
 & & \left|(f_n(\tau)+K)\hat{c}_n +
  \sum_{k=1}^{\infty} \frac{K}{(2k+1)!} \left((x+\hat{c})^{*(2k+1)}\right)_n\right|\\
 & = & |\left(\left(f_n(\tau)+K\right)\hat{c}_n
- \frac{K}{3!}\left(\hat{c}^{*3}\right)_n
+ \frac{K}{5!}\left(\hat{c}^{*5}\right)_n
- \frac{K}{7!}\left(\hat{c}^{*7}\right)_n\right) +\\
 & & \sum_{p=1}^{7}\sum_{j=-6l}^{6l} \gamma_{p, j} x^{*k}_{n-j} + R(\hat{c}+x) | \\
 & \leq & \left|\left(f_n(\tau)+K\right)\hat{c}_n
- \frac{K}{3!}\left(\hat{c}^{*3}\right)_n
+ \frac{K}{5!}\left(\hat{c}^{*5}\right)_n
- \frac{K}{7!}\left(\hat{c}^{*7}\right)_n\right| +\\
 & & + \left| \sum_{p=1}^{7}\sum_{j=-6l}^{6l} \gamma_{p, j} x^{*k}_{n-j} \right| + |R(\hat{c}+x)| \\
& \leq & \left|\left(f_n(\tau) - f_n(\tau_0)\right) \hat{c}_n \right| + \\
 & & + \left|\left(f_n(\tau_0)+K\right)\hat{c}_n
- \frac{K}{3!}\left(\hat{c}^{*3}\right)_n
+ \frac{K}{5!}\left(\hat{c}^{*5}\right)_n
- \frac{K}{7!}\left(\hat{c}^{*7}\right)_n\right| +\\
 & & + \sum_{p=1}^{7}\sum_{j=-6l}^{6l} \left|\gamma_{p, j}\right|
\frac{\beta_2(10\beta_2)^{p-1}}{(|n-j|+1)^2} \\
 & & + \sum_{k=4}^{22} \left| \left((x+\hat{c})^{*(2k+1)}\right)_n\right|
     + \left|\sum_{k=23}^{\infty} \left((x+\hat{c})^{*(2k+1)}\right)_n\right|
\end{array}
\end{equation*}

in the first term of RHS we will use the estimation:
\begin{equation}
\label{eq:zmianaf}
 \begin{array}{rcl}

\left|\left(f_n(\tau) - f_n(\tau_0)\right) \right| & \leq &
\max_{\tau \in [\underline{\tau}, \overline{\tau}]}
| f'(\tau)(\tau-\tau_0) | \leq \\
& \leq & \max_{\tau \in [\underline{\tau}, \overline{\tau}]}
| ine^{in\tau} - n^2\tau e^{in\tau}| \cdot |\Delta \tau| \leq \\
& \leq & \max_{\tau \in [\underline{\tau}, \overline{\tau}]}
(n + n^2\tau) |\Delta \tau|
   \end{array}
\end{equation}

The fact that it grows quickly with $n$ is irrelevant as for $n>5$ we have $\hat{c}_n=0$.

The
$\left|\left(f_n(\tau_0)+K\right)\hat{c}_n
- \frac{K}{3!}\left(\hat{c}^{*3}\right)_n
+ \frac{K}{5!}\left(\hat{c}^{*5}\right)_n
- \frac{K}{7!}\left(\hat{c}^{*7}\right)_n\right|$ can be directly computed
and is small as it is the numerical solution that is close to zero.

The $\sum_{p=1}^{7}\sum_{j=-6l}^{6l} \left|\gamma_{p, j}\right|
\frac{\beta_2(10\beta_2)^{p-1}}{(|n-j|+1)^2}$
can be directly computed and is small for $p>1$ as $10\beta_2$ is small.
If the terms for $p=1$ wouldn't be small we could try to diagonalize linear part. But there is no such need -- we check that the inequality holds.
We estimate $\left| \left((x+\hat{c})^{*(2k+1)}\right)_n\right|$ for $k \in \{4, .., 22\}$ from Lemma \ref{lem:srodek}. To estimate $\left|\sum_{k=23}^{\infty} \left((x+\hat{c})^{*(2k+1)}\right)_n\right|$ we use Lemma \ref{lem:ogon}. In
both estimates we have $\frac{1}{(2k+1)!}$ with $k$ large enough to make the result
small.

Thus we see that unless the $\gamma_{1, \cdot}$ are big, the
estimation of RHS should be small. The program computes the
estimation for LHS and RHS, compares them and for each $n$ in $0,
2, 3, 4, .., 224$ the inequality to hold. This ends the case for
$n \neq 1$.

\begin{table}
\caption{Estimates of LHS and RHS for small $n$}
\begin{center}
\begin{tabular}{|c|c|c|}
\hline
n & LHS & RHS \\
\hline
0 & $4.1792 \cdot 10^{-7}$ & $0.1687 \cdot 10^{-7}$ \\
2 & $0.4474 \cdot 10^{-7}$ & $0.0958 \cdot 10^{-7}$ \\
3 & $0.5080 \cdot 10^{-7}$ & $0.0827 \cdot 10^{-7}$ \\
4 & $0.4893 \cdot 10^{-7}$ & $0.0160 \cdot 10^{-7}$ \\
5 & $0.4537 \cdot 10^{-7}$ & $0.0118 \cdot 10^{-7}$ \\
6 & $0.4171 \cdot 10^{-7}$ & $0.0070 \cdot 10^{-7}$ \\
7 & $0.3834 \cdot 10^{-7}$ & $0.3465 \cdot 10^{-7}$ \\
8 & $0.3536 \cdot 10^{-7}$ & $0.0040 \cdot 10^{-7}$ \\
9 & $0.3274 \cdot 10^{-7}$ & $0.0036 \cdot 10^{-7}$ \\
\hline
\end{tabular}
\end{center}

\vspace{0.2cm}
{\footnotesize The estimates obtained by the program. The LHS is closest to the RHS for $n=7$ because we have artificially chosen $\hat{c}_7 = 0$. That makes the term
$\left|\left(f_n(\tau)+K\right)\hat{c}_n - \frac{K}{3!}\left(\hat{c}^{*3}\right)_n + \frac{K}{5!}\left(\hat{c}^{*5}\right)_n - \frac{K}{7!}\left(\hat{c}^{*7}\right)_n\right|$ for $n=7$ equal to approximately $0.3411 \cdot 10^{-7}$.}
\end{table}

{\bf 2.} We have two cases left: $\tau \in \{\underline{\tau}, \overline{\tau}\}$
and $|x_1| = \frac{\beta_2}{4}$. We have two variables left but only one
equation -- the equation for $n=1$. However this is a complex equation and the
variables are real so we will be able to show that in both cases $H(\tau, x) \neq 0$.

Let's denote:

\[ L_\tau := L_1(\tau)\hat{c}_1 \]
\[ L_x := \left(f_n(\tau_0)+K+\gamma_{1, 0}+\gamma_{1, 2} \right)x_1 \]
\[ \begin{array}{rcl}
    N & := & \left(f_n(\tau_0)+K\right)\hat{c}_n
- \frac{K}{3!}\left(\hat{c}^{*3}\right)_n
+ \frac{K}{5!}\left(\hat{c}^{*5}\right)_n
- \frac{K}{7!}\left(\hat{c}^{*7}\right)_n + \\
      & & + \sum_{-6l \leq j \leq 6l, j \neq 0,2} \gamma_{p, j} x^{*k}_{n-j} +
        \sum_{p=2}^{7}\sum_{j=-6l}^{6l} \gamma_{p, j} x^{*k}_{n-j} + \\
      & & + R(c) + r_n(\tau)\hat{c}_n + (f(\tau)-f(\tau_0))x_n
   \end{array}
\]

The $L_\tau, L_x, N$ depend on $\tau, x$ but to make the notation short we won't write it.  They are chosen such that $G(\tau, x)_1 = L_\tau + L_x$ and $F(\tau, x)_1 = L_\tau + L_x + N$ and $H(h, \tau, x)_1 = L_\tau + L_x + hN$. If $\tau \in \{\underline{\tau},
\overline{\tau}\}$ then we show that:

\[
\left| L_\tau \right| > \left|L_x\right| + \left|N\right|
\]

The $\left| L_\tau \right|$ is estimated as:

\[ |L_1(\tau)\hat{c}_1| = \left. \begin{array}{l} \\ |in-n^2|\cdot \Delta\tau \cdot |\hat{c}_n| \\~ \end{array} \right|_{n=1} =
\sqrt{2}\Delta\tau \cdot |\hat{c}_1| \]

We estimate the terms in $|N|$ as for $n \neq 1$ with the exception
for the $\left|(f(\tau)-f(\tau_0))x_n\right|$ and $\left|r_n(\tau)\hat{c}_n\right|$ that were not present in the previous case. The the first one we use the inequality (\ref{eq:zmianaf}). For the second we use the estimation:

\[
\begin{array}{rcl}
|r_n(\tau)\hat{c}_n| & \leq &
\sup_{\tau \in [\underline{\tau};\overline{\tau}]}
\left|\frac{1}{2}f_1''(\tau)(\tau-\tau_0)^2\right|\cdot |\hat{c}_n| \\
&=&\frac{\left|\hat{c}_n\right|}{2} \sup_{\tau \in [\underline{\tau};\overline{\tau}]}
\left| -\left(2n^2 + in^3\tau\right) e^{in\tau}\right| \\
&=&\frac{\left|\hat{c}_n\right|}{2} \sup_{\tau \in [\underline{\tau};\overline{\tau}]}
\left| 2 + i\tau \right|
\end{array}
\]

Having all these estimates our program checks that the inequality
is satisfied.

In the case $|x_1| = \frac{\beta_2}{4}$ the inequality $|L_x| > |L_\tau| + |N|$
is of course false. To prove that there is no zero will need to use the fact that
$\tau \in \setR$ and $x_1 \in \setR$. On the figure 1 we have a sketch how
the sets $L_\tau$ and $L_x + N$ looks like -- we see that they shouldn't intersect.

\begin{figure}
\psfrag{Possible values of Ltau}{$L_\tau$}
\psfrag{Possibe value of LxN}{$L_x+N$}
\includegraphics[width=10cm]{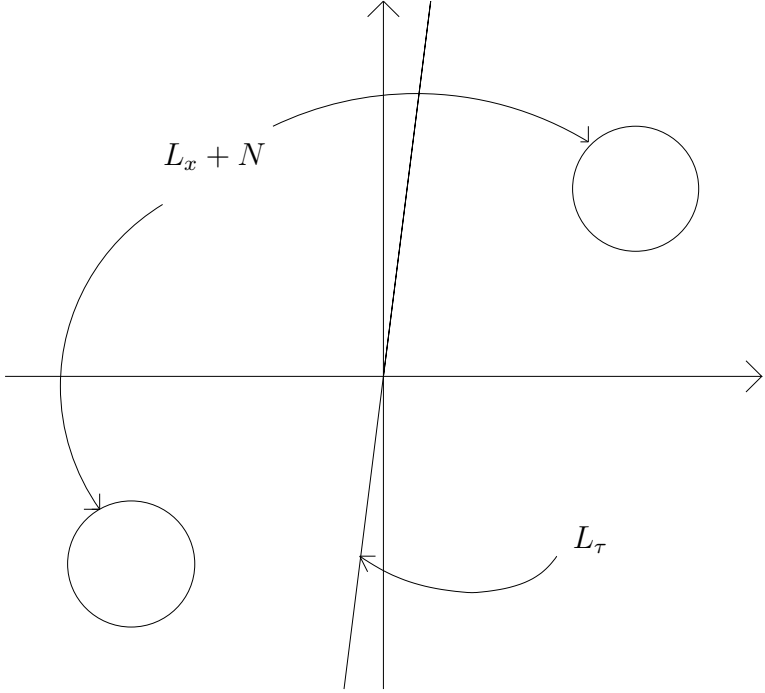}
\caption{The possible values of $L_\tau$ and $L_x+N$}
\end{figure}

Formally we will show that $\{\tan\arg(L_\tau):
(\tau, x) \in [\underline{\tau}, \overline{\tau}] \times X_3,
x_1 = \pm \frac{\beta_2}{4} \} \cap \{\tan\arg(L_x + N):
(\tau, x) \in [\underline{\tau}, \overline{\tau}] \times X_3,
x_1 = \pm \frac{\beta_2}{4} \} = \emptyset$ (where $\arg$ is the
complex number argument). Of course this implies that the sum of
these terms is different from zero.

The $\arg(L_\tau)$ is easy to compute as this is $\tan\arg\left(\left(ie^{i\tau_0}-e^{i\tau_0}\right)c_1\right)$
and $\tan\arg z$ for a complex number $z$ can be computed as
$\frac{\Im z}{\Re z}$. To estimate the other set we will use the estimation
for $N$ from the previous point -- let's assume that $|N| \leq \lambda$
for some $\lambda > 0$. We have that $L_x + N \in
\pm \left(f_n(\tau)+K+\gamma_{1, 0}+\gamma_{1, 2}\right)\frac{\beta_2}{4}
+[-\lambda; \lambda] + [-\lambda; \lambda]i$ and using the interval
arithmetics we can find $\tan \arg(L_x + N)$. The program checks that these two
sets are disjoint and this ends the proof of the theorem.
\qed

To finish the proof of the existence of the orbit let's define a
second homotopy.

\[ H^{L}(h, \tau, x) := \left\{
\begin{array}{l}
\left(f_n(\tau_0+h(\tau-\tau_0)) + K\right)x_n\mbox{ for }n \neq \pm 1\\
\left(f_n(\tau_0) + K + \gamma_{1, 0} + \gamma_{1,2n}\right)x_n + L(\tau)c_1\mbox{ for }n = \pm 1
\end{array}\right.\]

It deforms $G$ into:

\[ G^{L}(h, \tau, x) := \left\{
\begin{array}{l}
\left(f_n(\tau_0) + K\right)x_n\mbox{ for }n \neq \pm 1\\
\left(f_n(\tau_0) + K + \gamma_{1, 0} + \gamma_{1,2n}\right)x_n + L(\tau)c_1\mbox{ for }n = \pm 1
\end{array}\right.\]

\begin{lemma}
If $(\tau, x) \in [\underline{\tau}, \overline{\tau}]\times X_3$ such that $\tau \in \{\underline{\tau}, \overline{\tau}\}$
or $\exists n: \left|x_n\right| = \frac{\beta_2}{(|n|+1)^2}$ then $H^L(h, x, \tau) \neq 0$
\end{lemma}

{\bf Proof:} Let $h, x, \tau$ be such that $H^L(h, x, \tau) = 0$.

Let $n \neq \pm 1$. The $\tau_0+h(\tau-\tau_0) \in [\underline{\tau}; \overline{\tau}]$
thus $|f_n(\tau_0+h(\tau-\tau_0)) + K|$ can be estimated as
in equation (\ref{eq:homotoplhs}). For each such value we have proven that it is
strictly greater than an RHS$\geq 0$.
Hence $|f_n(\tau_0+h(\tau-\tau_0)) + K| > 0$ and $H^L(h, \tau, x) = 0 \Rightarrow x_n = 0$.

But if $x_n = 0$ for each $n \neq \pm 1$ then $H^L(h, \tau, x) = G(\tau, x) = H(0, \tau, x)$ and $H$ has no zeros on the boundary. \qed

\begin{observation}
The Galerkin projection of $H$: i.e. $H_l(\tau, y_0, \dots, y_l)
:= P_l H(\tau, Q_l(y_0, \dots, y_l))$ (for $(y_0, ..., y_n) \in
\setR \times \setR \times \setC^{n-1}$) is a homotopy from
$F^l(\tau, c+\cdot)$ to the projection of $G$. Analogically the
projection of $H^L$ is a homotopy from $G$ to $G^L$. There are no
zeros on the boundaries for these homotopies.
\end{observation}

\begin{observation}
The degree $\deg(F_l, [\underline{\tau}, \overline{\tau}] \times P_l(X_3), 0)$ is well defined.
\end{observation}

{\bf Proof:}
The homotopy $H_l$ for $h=1$ have no zeros on the boundary and this is
$F_l$. \qed

\begin{lemma}
\label{lem:stopien}
$\deg(F_l, [\underline{\tau}, \overline{\tau}] \times P_l(X_3), 0) \neq 0$
\end{lemma}

{\bf Proof:} We know that the degree of $F_l$ is equal to the
degree of the Galerkin projection of $G^L$ -- let's denote it by
$G^L_l$. The $G^L_l$ is a linear function. If the determinant of
the differential was zero then there would be a zero on the
boundary of any neighborhood of $(\tau_0, 0)$ so the determinant
is non-zero. That means that the degree is $\pm 1$ thus it is
non-zero. \qed

\begin{theorem}
The equation (\ref{eq:rownanie}) has a
periodic solution for some $\tau \in [\underline{\tau},
\overline{\tau}]$ whose Fourier coefficients are in the set $X_3$.
\end{theorem}

{\bf Proof: } From the Lemma \ref{lem:stopien} we know that the
local Brower degree is non-zero thus each $F_l$ has a zero in
$X_3$. From Lemma \ref{lem:granica} we have that $F$ has a zero in
$X_3$. From Theorem \ref{tw:fourier} we obtain that in $X_3$ there
is a solution of the equation. \qed

\section{Conclusions}

In this paper I was able to rigorously prove the existance of a periodic orbit for $K=1.6$. I was not able to show that the period is 4 even if the numerical simulations suggests that -- in the proof we need the $\tau$ as a variable for the image and the domain to have the same dimension.

The value $K = 1.6$ has been chosen because it's easiest to find the $\hat{c}_n$ values for an attracting orbit. However using the Newton method it should be possible to find an approximated orbit which is not attracting (from the numerical simulations it seems to happen for $K > 5.11$). For larger $K$ the values of $\hat{c}_n$ may not decay as fast so we may need to diagonalize the first coefficients of $F_l$ as mentioned in the proof.

Of course $\sin$ is periodic so there exist infintely many such orbit that differs by $2k\pi$ ($k \in \setZ$). As mentioned for $K > 5.11$ is seems that these periodic orbits stops to be attracting and orbits from numerical simulations jump $\pm 2\pi$ from one periodic orbit to another. This suggests heteroclinic connections and chaos but proving it would require some new ideas.

\section{Acknowledgements}

I would like to thank Piotr Zgliczy\'nski for suggesting the topic and helpful discussions.


\begin{thebibliography}{99}
\bibitem[Cesari 64]{cesari} L. Cesari, Functional analysis and Galerkin methods, {\em Michigan Math. J.} {\bf 11} (1964) 383-414
%
\bibitem[Wischert et al 94]{fizycy} W. Wischert, W. Wunderlin, A. Pelster, M. Olivier, and J. Groslambert
{\em Delay-induced instabilities in nonlinear feedback systems}, Physical Review E {\bf 49} 203-219 (1994)
%
\bibitem[Smoller 83]{S} J. Smoller, {\em Shock Waves and Reaction-Diffusion
Equations} Springer, 1983.
%
\bibitem[Zgliczy\'n{}ski 01]{zgliczyn2001} P. Zgliczy\'nski, K. Mischaikow
{\em Rigorous Numerics for Partial Differential Equations: The Kuramoto-Sivashinsky Equation}, Found. Comput. Math {\bf 1} 255-288 (2001)
%
\bibitem[Zgliczy\'n{}ski 04]{zgliczyn} P. Zgliczy\'nski
{\em Rigorous Numerics for Dissipative Partial Differential Equations II. Periodic Orbit for the Kuramoto-Sivashinsky Equation}, Found. Comput. Math 157-195 (2004)
\end{thebibliography}
\end{document}